\documentclass[a4paper]{article}
\usepackage[pdftex]{graphicx}
\usepackage{amsmath}
\usepackage{amssymb}
\usepackage{amsthm}
\usepackage[utf8]{inputenc}
\usepackage[T1]{fontenc}
\usepackage{hyperref}
\usepackage{smartref}
\usepackage{todonotes}
\usepackage{fullpage}
\usepackage{epstopdf}
\usepackage{subcaption}
\usepackage[inline]{enumitem}

\newtheorem{lemma}{Lemma}
\newtheorem{theorem}{Theorem}

\usepackage[toc,page]{appendix}
\usepackage{algorithm}
\usepackage{algorithmic}

\usepackage{color}

\renewcommand\vec{\boldsymbol}

\newcount\colveccount
\newcommand*\colvec[1]{
        \global\colveccount#1
        \begin{pmatrix}
        \colvecnext
}
\def\colvecnext#1{
        #1
        \global\advance\colveccount-1
        \ifnum\colveccount>0
                \\
                \expandafter\colvecnext
        \else
                \end{pmatrix}
        \fi
}

\begin{document}
%\pagecolor{black}
%\color{white}

\title{On generic $3$-rigidity of graphs}
\author{Tam\'as Baranyai}

\maketitle

\abstract{We give a necessary condition of generic $3$ -rigidity of graphs relying on partitioning the edges into $3$ subsets; such that each subset-pair gives a generically $2$-rigid graph, either by themselves or after an appropriate edge-deletion. Notably, as pointed out by Dewar and Gallet, the condition is still not sufficient.}

\section{Introduction}
Let $G(V,E)$ be a graph such that there are no double edges. We will consider pin jointed bar and joint frameworks (trusses) having $G(V,E)$ as a defining graph, as follows. Let $\{ \vec{x}_1, \ \dots, \ \vec{x}_d \}$ be a base in $\mathbb{R}^d$ and $\vec{p}:$ $V \to \mathbb{R}^d$ be a map that assigns to each vertex a point in $d$ dimensional space in a generic way. We assume $|V|>d$ as this is the mechanically interesting case. Let there be a ball joint placed at each point, that will connect the bars of the structure meeting at these points. For each edge $e \in E$ consider a bar to run between the joints corresponding to the vertex pair of $e$. Graph $G(V,E)$ is called generically $d$-rigid if the bar and joint framework is rigid in $d$-dimensions. Rigidity depends on $G$ as well as the map $\vec{p}$ but since we specified the points to be in generic position it becomes dependant only on $G(V,E)$. We drop the term "generic" for brevity. Graph $G(V,E)$ is called minimally $d$-rigid if it is $d$-rigid and has the minimal number of edges to do so.

There is a large amount of literature on the rigidity of this type of structure, for the latest survey-like paper we refer the reader to the work of Cruickshank et. al. \cite{cruickshank2025rigidity}. A combinatorical characterization of rigidity has been found for $d=1$ and $d=2$ \cite{Pollaczek-Geiringer} and as far as we know it is a long-standing open problem for $d \geq 3$. 

The work intended to generalize the result of Lovász and Yemini:

\begin{theorem}[Lovász and Yemini \cite{Lovasz-Yemini}]
Graph $G(V,E)$ is minimally $2$-rigid if and only if after doubling any edge the resulting graph may be decomposed into two edge disjoint spanning trees.
\end{theorem}

As the three dimensional analogue of this is necessary but not sufficient, we looked for extra conditions on the decomposition. We found planar rigidity of pairs of edge-sets, either by themselves or after appropriate edge-deletions. While we originally thought this to be a necessary and sufficient condition, it turned out to be necessary only. As such we archive here the part of the original proof which aims for the necessity direction. To see how it is not sufficient we refer the reader to the work of Dewar \cite{dewar2026counterexamplebaranyaiscombinatorialcharacterisation}, which also gives an alternative proof.

\section{Self-stresses of frameworks}    

First we investigate the possible self-stresses of the framework, which are non-trivial internal force distributions satisfying static equilibrium in the absence of external loads; since $G$ is minimally $d$-rigid exactly if the truss has $d |V| -  {d+1 \choose 2}$ bars and no self-stresses.

If $\vec{B}$ is a signed incidence matrix of the graph, in the absence of external loads the equilibrium of force components in the $\vec{x}_i$ direction may be described as $\vec{B}\vec{s}^i=\vec{0}$ where vector $\vec{s}^i$ collects the components of the bar forces corresponding to each edge. The equilibrium of force components may be used to re-parametrize the internal force-distribution of the truss the following way: Delete the row corresponding to vertex $v_1$ from $\vec{B}$ and partition the remaining matrix into $[\vec{F} | \vec{R}]$ where $\vec{F}$ is invertible and corresponds to an arbitrary spanning tree $F$, while matrix $\vec{R}$ contains the edges of the remainder $R=V \setminus F$. The static equilibrium becomes $\vec{F}\vec{f}^i+\vec{R}\vec{r}^i=\vec{0} \Leftrightarrow \vec{f}^i=-\vec{F}^{-1}\vec{R}\vec{r}^i$. To have this formulation describe the truss we have to enforce on the bar forces the directional constraints of the bars. This may be done with the help of diagonal matrices $\vec{D}_f^{i,j}$ and $\vec{D}_r^{i,j}$ to enforce $\vec{f}^j=\vec{D}_f^{i,j}\vec{f}^i$ and $\vec{r}^j=\vec{D}_r^{i,j}\vec{r}^i$ ($i,j \in \{1, \dots d \}$ and $i\neq j$). There are two types of redundancies present in this description. The first is due to the fact that the diagonal matrices over-parametrize the geometry as $\vec{r}^j=\vec{D}_r^{i,j}\vec{r}^i$ and $\vec{r}^k=\vec{D}_r^{j,k}\vec{r}^j$ together imply $\vec{r}^k=\vec{D}_r^{i,k}\vec{r}^i$ and a similar argument holds for the spanning-tree edges. This will be handled by restricting the indices: we will choose the set with $i=1$ and $j$ running from $2$ to $d$. As such, the number of parameter-vectors $\vec{r}^i$ is reduced from $d$ to $1$ and the equations they are subject to are given by the geometry of the edges of $F$ as
\begin{align}
\left[\vec{D}_f^{1,j} \vec{F}^{-1}\vec{R}- \vec{F}^{-1}\vec{R} \vec{D}_r^{1,j} \right]\vec{r}^1=\vec{0} \label{eq:1}
\end{align} 
for $j \in \{ 2, \dots d \}$. The second type of redundancy comes from the fact that the directional constraints of the bar forces are consequences of moment-equilibrium equations; implying one redundancy for each coordinate plane. Formally we have the following lemma:

\begin{lemma}\label{L:0}
For any $e\in E$ and any coordinate plane $[\vec{x}_k,\vec{x}_l ]$ if  $\vec{B}\vec{s}^k=\vec{0}$ and $\vec{B}\vec{s}^l=\vec{0}$ are known to hold and the directional constraints in the plane are known to be satisfied for all edges except for $e$ then it holds for $e$ as well.
\end{lemma} 

\begin{proof}
Let us collect the coordinate-differences of the endpoints of the bars into vectors $\vec{g}^k \in \mathbb{R}^{|E|}$ for each $k$ as follows. Using the notation that $p_k^j$ is the $k$-th coordinate of $\vec{p}(v_j)$, the entry corresponding to $e_{i,j}$ in $\vec{g}^k$ is $p_k^j-p_k^i$. Let us partition  $\vec{g}^k$ into $\vec{g}^k_f$ and $\vec{g}^k_r$ according to the spanning tree and the remainder.  We have
\begin{align}
\vec{g}^k_f \left[\vec{D}_f^{k,l} \vec{F}^{-1}\vec{R}- \vec{F}^{-1}\vec{R} \vec{D}_r^{k,l} \right]= \vec{g}^l_r-\vec{g}^l_r =\vec{0}. 
\end{align}
Since no entry of $\vec{g}^k_f$ is zero any one row of $\left[\vec{D}_f^{k,l} \vec{F}^{-1}\vec{R}- \vec{F}^{-1}\vec{R} \vec{D}_r^{k,l} \right]$ may be expressed using all the others. 
\end{proof}

\section{$3$-rigidity of graphs}
In the following we focus on the case of $d=3$, denote the deletion of the top row of a matrix $\vec{A}$ by $\overline{\vec{A}}$, and denote edge-contraction with $/e$. With some abuse of notation we will use this notation to contract edges that have been removed from the graph, which is to be interpreted as vertex identification. We need two lemmas before the main result. 
\begin{lemma}\label{L:1}
If $G(V,E)$ is minimally $3$-rigid then for any spanning tree $F \subset E$ and any edge $e\in F$ there is a partition of $V \setminus F$ into $R_1$ and $R_2$ such that $G(V,F \cup R_1)$ and $G( V,F \cup R_2)/ e $ are minimally $2$-rigid.
\end{lemma}

\begin{proof}
Without the loss of generality we may order the rows of $\vec{F}^{-1}\vec{R}$ such that the top row corresponds to $e$, and assume a coordinate system where the bar corresponding to $e$ is parallel to $\vec{x}_1$. According to Lemma \ref{L:0} we may delete the first rows of Equation \ref{eq:1} for indices $j=2$ and $j=3$ without losing an independent equation. Now, consider the projection of the $3$-dimensional problem to the $[\vec{x}_2,\vec{x}_3 ]$ plane. Since the directional constraints of edge $e$ play no role in this particular planar projection of the structure, we may conclude using Lemma \ref{L:0} that if all the directional constraints are known to hold in the $[\vec{x}_2,\vec{x}_3 ]$ plane except for that of $e$ and an additional single edge, then the constraint holds for the additional single edge. (A linear algebraic reasoning would be to note that the first entry of $\vec{g}^2_f$ is zero, yet $\vec{g}^2_f \left[\vec{D}_f^{2,3} \vec{F}^{-1}\vec{R}- \vec{F}^{-1}\vec{R} \vec{D}_r^{2,3} \right]=\vec{0}$.) Since we have chosen to express all spatial constrains in the $[\vec{x}_1,\vec{x}_2 ]$ and $[\vec{x}_1,\vec{x}_3 ]$ planes, we have to delete any one row from the two coefficient matrices coming from Equation \ref{eq:1} additionally to the rows corresponding to edge $e$. we conclude that $\vec{r}^1$ may be determined using the system of linear equations
\begin{align}
\begin{bmatrix}
\overline{\vec{D}_f^{1,2} \vec{F}^{-1}\vec{R}- \vec{F}^{-1}\vec{R} \vec{D}_r^{1,2}}\\
\overline{\overline{\vec{D}_f^{1,3} \vec{F}^{-1}\vec{R}- \vec{F}^{-1}\vec{R} \vec{D}_r^{1,3}}} 
\end{bmatrix} \vec{r}^1=\vec{0} \label{eq:2}
\end{align}
the coefficient-matrix of which is invertible if and only if $G(V,E)$ is minimally $3$-rigid. Looking at the Laplace-expansion of the non-zero determinant along the first $n-2$ rows shows that if $G$ is $3$-rigid there must be a partition $\vec{R}=[\vec{R}_1|\vec{R}_2]$ which extends to the diagonal matrices and induces a partition of edges $R=R_1 \cup R_2$; such that $\overline{\vec{D}_f^{1,2} \vec{F}^{-1}\vec{R}_1- \vec{F}^{-1}\vec{R}_1 \vec{D}_{r1}^{1,2}}$ and $\overline{\overline{\vec{D}_f^{1,3} \vec{F}^{-1}\vec{R}_2- \vec{F}^{-1}\vec{R}_2 \vec{D}_{r2}^{1,3}}} $ are invertible. Considering how these equations would look like in case of $d=2$ show us that the first invertibility condition means graph $G(V,F \cup R_1)$ is minimally $2$-rigid. To interpret the second invertibility condition we observe that rows of $[\vec{F}^{-1}\vec{R]}$ correspond to edges of $F$ while columns of it to cycles of $G$. Deleting a row from  $[\vec{F}^{-1}\vec{R]}$ is an operation that reduces the number of spanning tree edges by one without altering the number of cycles or what other edges the cycles traverse. This may be achieved by contracting the spanning tree-edge corresponding to the deleted row. Since the diagonal matrices or the subtraction in Equation \eqref{eq:2} does not alter the structure of the matrices involved, we may interpret the invertibility of $\overline{\overline{\vec{D}_f^{1,3} \vec{F}^{-1}\vec{R}_2- \vec{F}^{-1}\vec{R}_2 \vec{D}_{r2}^{1,3}}} $ as the graph $G( V,F \cup R_2)/ e$ being minimally $2$-rigid. Since $G(V,E)$ has no double-edges the contraction of $e$ will not create loop-edges and this interpretation makes sense. 
\end{proof}

In order to proceed we need to rely on rigidity matrices. Let $\vec{M}(G,d)$ denote the the $d$-dimensional rigidity matrix of graph $G(V,E)$ such that rows correspond to edges of $G$, columns are partitioned into blocks $\vec{C}_i$ corresponding to the directions  $\vec{x}_i$. In each block $\vec{C}_i$ the columns correspond to the vertices of $G$ and in each row corresponding to edge $e_{j,k}$ there are two entries: $p_i^k-p_i^j$ in the $j$-th column and $p_i^j-p_i^k$ in the $k$-th column ($p_i^j$ is the $i$-th coordinate of $\vec{p}(v_j)$ ). Graph $G$ is minimally $d$-rigid exactly if $|E|= d |V| -  {d+1 \choose 2}$ and $\vec{M}(G,d)$ has full (row) rank. In what follows we will omit $G$ and/or $d$ from the notation of the rigidity matrix if the context does not require such distinctions. In order to use determinants of square matrices we delete the first $i$ columns of block $\vec{C}_i$ for all $i$ and denote the result with $\vec{N}(G,d)$. It is easy to see for $d=2$ and for $d=3$ that if the points $\vec{p}(v_i)$ are in general position with respect to the coordinate system then this is a good way to support a rigid body and $\vec{N}(G,2)$ and $\vec{N}(G,3)$ are invertible exactly if $G$ is minimally $2$ and $3$-rigid respectively. Although we specified the map $\vec{M} \mapsto \vec{N}$ using the labelling of vertices its properties hold for any labelling. We may now continue with the following Lemma:

\begin{lemma}\label{L:2}
If $G(V,E)$ is minimally $3$-rigid then for each and any $e \in E$ there is a spanning tree $F \subset E$  such that $e\in F$ and $G( V,E \setminus F)/ e $ is minimally $2$-rigid.
\end{lemma}

\begin{proof}
Without the loss of generality we may label the vertices such that the specified $e \in E$ runs between vertices $v_1$ and $v_2$ and assume a coordinate system where the bar corresponding to $e$ is parallel to $\vec{x}_1$. Since $G(V,E)$ is minimally $3$-rigid the Laplace expansion of the determinant $|\vec{N}(G,3)|$ along the first $n-1$ columns proves the existence of two invertible blocks, the $(n-1) \times (n-1)$ sized $\vec{N}_{11}$ and its complement $\vec{N}_c$. The invertibility of $\vec{N}_{11}$  means its rows correspond to spanning tree $F$ in $G$. Due to the assumption on the coordinate system the complement $\vec{N}_c$ is of shape $\vec{N}(G( V,E \setminus F)/ e,2)$ thus the remaining graph is minimally $2$-rigid.  
\end{proof}

Our main result is the following Theorem:

\begin{theorem}
If graph $G(V,E)$ is minimally $3$-rigid then for any $e \in E$ it may be decomposed into graphs $G(V,S_1),G(V,S_2),G(V,S_3)$ such that $|S_i|=|V|-i$ and the graphs $G(V,S_1 \cup S_2)$,  $G(V,S_1 \cup S_3) / e$ and $G(V,S_2 \cup S_3) / e$ are minimally $2$-rigid. 
\end{theorem}

\begin{proof}
To see how rigidity implies the decomposition one may use Lemma \ref{L:2} to get a spanning tree which will be $S_1$ then use Lemma \ref{L:1} to decompose the remainder.
\end{proof}

\section{Acknowledgements}
We are grateful to Sean Dewar and Matteo Gallet for pointing out that the original claim is only partially correct.

\bibliographystyle{unsrt} 
\bibliography{rigidbib}

\begin{thebibliography}{1}

\bibitem{cruickshank2025rigidity}
James Cruickshank, Bill Jackson, Tibor Jord{\'a}n, and Shin-ichi Tanigawa.
\newblock Rigidity of graphs and frameworks: A matroid theoretic approach.
\newblock {\em arXiv preprint arXiv:2508.11636}, 2025.

\bibitem{Pollaczek-Geiringer}
H.~Pollaczek-Geiringer.
\newblock Über die gliederung ebener fachwerke.
\newblock {\em ZAMM - Journal of Applied Mathematics and Mechanics /
  Zeitschrift für Angewandte Mathematik und Mechanik}, 7(1):58--72, 1927.

\bibitem{Lovasz-Yemini}
L.~Lov\'{a}sz and Y.~Yemini.
\newblock On generic rigidity in the plane.
\newblock {\em SIAM Journal on Algebraic Discrete Methods}, 3(1):91--98, 1982.

\bibitem{dewar2026counterexamplebaranyaiscombinatorialcharacterisation}
Sean Dewar.
\newblock A counter-example to Baranyai's combinatorial characterisation for
  3-rigidity.
\newblock {\em arXiv preprint arXiv:2601.19460}, 2026.

\end{thebibliography}

\end{document}